\magnification=1100
\baselineskip=14truept
\voffset=.75in  
\hoffset=1truein
\hsize=4.5truein
\vsize=7.75truein
\parindent=.166666in
\pretolerance=500 \tolerance=1000 \brokenpenalty=5000


\def\anote#1#2#3{\smash{\kern#1in{\raise#2in\hbox{#3}}}%
  \nointerlineskip}     
\def\note#1{%
  \hfuzz=50pt%
  \vadjust{%
    \setbox1=\vtop{%
      \hsize 3cm\parindent=0pt\eightrm\baselineskip=9pt%
      \rightskip=4mm plus 4mm\raggedright#1%
      }%
    \hbox{\kern-4cm\smash{\box1}\hfil\par}%
    }%
  \hfuzz=0pt
  }
\def\note#1{\relax}

\newcount\equanumber
\equanumber=0
\newcount\sectionnumber
\sectionnumber=0
\newcount\subsectionnumber
\subsectionnumber=0
\newcount\snumber  
\snumber=0

\def\section#1{%
  \subsectionnumber=0%
  \snumber=0%
  \equanumber=0%
  \advance\sectionnumber by 1%
  \noindent{\bf \the\sectionnumber .~#1~}%
}%
\def\subsection#1{%
  \advance\subsectionnumber by 1%
  \snumber=0%
  \equanumber=0%
  \noindent{\bf \the\sectionnumber .\the\subsectionnumber .~#1~}%
}%
\def\prevs{\the\sectionnumber .\the\snumber }
\long\def\Definition#1{\noindent{\bf Definition.}{\it#1}}

\long\def\Corollary#1{%
  \global\advance\snumber by 1%
  \bigskip
  \noindent{\bf Corollary~\prevs .}%
  \quad{\it#1}%
}%
\long\def\Lemma#1{%
  \global\advance\snumber by 1%
  \bigskip
  \noindent{\bf Lemma~\prevs .}%
  \quad{\it#1}%
}%
\def\Proof{\noindent{\bf Proof.~}}
\long\def\Proposition#1{%
  \advance\snumber by 1%
  \bigskip
  \noindent{\bf Proposition~\prevs .}%
  \quad{\it#1}%
}%
\long\def\Remark#1{%
  \bigskip
  \noindent{\bf Remark.~}#1%
}%
\long\def\Theorem#1{%
  \advance\snumber by 1%
  \bigskip
  \noindent{\bf Theorem~\prevs .}%
  \quad{\it#1}%
}%
\def\ifundefined#1{\expandafter\ifx\csname#1\endcsname\relax}
\def\labeldef#1{\global\expandafter\edef\csname#1\endcsname{\prevs}}
\def\labelref#1{\expandafter\csname#1\endcsname}
\def\label#1{\ifundefined{#1}\labeldef{#1}\note{$<$#1$>$}\else\labelref{#1}\fi}

\def\preveq{(\the\sectionnumber .\the\equanumber)}
\def\neq{\global\advance\equanumber by 1\eqno{\preveq}}

\def\ifundefined#1{\expandafter\ifx\csname#1\endcsname\relax}

\def\equadef#1{\global\advance\equanumber by 1%
  \global\expandafter\edef\csname#1\endcsname{\preveq}%
  \preveq}

\def\equaref#1{\expandafter\csname#1\endcsname}
 
\def\equa#1{%
  \ifundefined{#1}%
    \equadef{#1}%
  \else\equaref{#1}\fi}

%
%
\font\eightrm=cmr8%
\font\sixrm=cmr6%
\font\fourrm=cmr4%
\font\eightbf=cmb8%
\font\elevenbb=msbm10%
\font\eigthbb=msbm8%
\font\sixbb=msbm6%
\newfam\bbfam%
\textfont\bbfam=\elevenbb%
\scriptfont\bbfam=\eigthbb%
\scriptscriptfont\bbfam=\sixbb%
\font\tencmssi=cmssi10%
\font\eightcmssi=cmssi7%
\font\sixcmssi=cmssi6%
\newfam\ssfam%
\textfont\ssfam=\tencmssi%
\scriptfont\ssfam=\eightcmssi%
\scriptscriptfont\ssfam=\sixcmssi%
\def\ssi{\fam\ssfam\tencmssi}%

\font\tenmsam=msam10%
\font\sevenmsam=msam7%
\font\sixmsam=msam6%

\def\bb{\fam\bbfam\elevenbb}%

\def\hexdigit#1{\ifnum#1<10 \number#1\else%
  \ifnum#1=10 A\else\ifnum#1=11 B\else\ifnum#1=12 C\else%
  \ifnum#1=13 D\else\ifnum#1=14 E\else\ifnum#1=15 F\fi%
  \fi\fi\fi\fi\fi\fi}
\newfam\msamfam%
\textfont\msamfam=\tenmsam%
\scriptfont\msamfam=\sevenmsam%
\scriptscriptfont\msamfam=\sixmsam%
\mathchardef\leq"3\hexdigit\msamfam 36%
\mathchardef\geq"3\hexdigit\msamfam 3E%

\def\convM{\buildrel{\raise-1pt\hbox{\sixrm M}}\over *}
\def\convMss{\setbox1=\hbox{\fourrm M}\ht1=0pt\dp1=0pt%
  \buildrel{\raise -2pt\box1}\over *}
\def\<{\langle}
\def\>{\rangle}

\def\d{\,{\rm d}}

\def\D{{\rm D}}
\long\def\DoNotPrint#1{\relax}

\def\Id{{\rm Id}}

\def\limt{\lim_{t\to\infty}}

\def\L{{\ssi L}}
\def\M{{\ssi M}}

\def\oF{\overline F{}}

\def\oG{\overline G{}}

\def\oK{\overline K{}}

\def\qed{~~~{\vrule height .9ex width .8ex depth -.1ex}}

\def\ss{\scriptstyle}
\def\T{\hbox{\tencmssi T}}

\def\HH{{\bb H\kern .5pt}}

\def\NN{{\bb N}\kern .5pt}

\def\RR{{\bb R}}

%
%
%
%
\def\uncatcodespecials 
    {\def\do##1{\catcode`##1=12}\dospecials}%
{\catcode`\^^I=\active \gdef^^I{\ \ \ \ }
 \catcode`\`=\active\gdef`{\relax\lq}}
\def\setupverbatim 
    {\parindent=0pt\tt %
     \spaceskip=0pt \xspaceskip=0pt 
     \catcode`\^^I=\active %
     \catcode`\`=\active %
     \def\par{\leavevmode\endgraf}
     \obeylines \uncatcodespecials \obeyspaces %
     }%
{\obeyspaces \global\let =\ }
%
%
%
%
%
%

\centerline{\bf ASYMPTOTIC EXPANSIONS FOR DISTRIBUTIONS}
\centerline{\bf OF COMPOUND SUMS OF}
\centerline{\bf LIGHT SUBEXPONENTIAL RANDOM VARIABLES}

\bigskip
 
\centerline{Ph.\ Barbe${}^{(1)}$, W.P.\ McCormick${}^{(2)}$ and 
C.\ Zhang${}^{(2)}$}
\centerline{${}^{(1)}$CNRS, France, and ${}^{(2)}$University of Georgia}
 
{\narrower
\baselineskip=9pt\parindent=0pt\eightrm

\bigskip

{\eightbf Abstract.} We derive an asymptotic expansion for the distribution
of a compound sum of independent random variables, all having the
same light-tailed subexponential distribution. The examples
of a Poisson and geometric number of summands serve as an illustration of
the main result. Complete calculations are done for a Weibull distribution,
with which we derive, as examples and without any difficulties, 7 terms 
expansions.

\bigskip


\bigskip
 
\noindent{\eightbf Keywords:} asymptotic expansion,
convolution, tail area approximation, regular variation, subexponential
distributions.

}

\bigskip\bigskip

\section{Introduction.}
In this
paper, we construct asymptotic expansions for the tail area $\oG$ of
a compound sum, when the summands belong to a class of light-tailed
subexponential distributions. To be more precise, let $X_i$, $i\geq 1$, 
be a sequence of independent random variables,
all having the same distribution $F$. For any positive integer $n$ the
partial sums $S_n=X_1+\cdots + X_n$ have distribution the $n$-fold
convolution $F^{\star n}$. We set $S_0=0$ and therefore $F^{\star 0}$
is defined as the distribution of the point mass at the origin. Let
$N$ be a nonnegative integer-valued random variable, independent of
the $X_i$'s. We consider the distribution $G$ of the compound sum
$S_N$, that is $EF^{\star N}$. Its tail area is $\oG=1-G$.
First order asymptotic results for $\oG$ have been obtained by
Embrechts, Goldie and Veraverbeke (1979),
Cline (1987), Embrechts (1985) and Gr\"ubel (1987).  A second order
formula may be found in  Omey and Willekens (1987).

Compound sums or
subordinated distributions arise as distribution of interest in
several stochastic models. In insurance risk theory, it models the
total claim amount. For a discussion of issues related to random sums
and insurance risk, we refer to Embrechts, Kl\"uppelberg and Mikosch
(1997), Asmussen (1997), Goldie and Kl\"uppelberg (1998). Compound
sums also appear in queueing theory, in connection with the stationary
distribution of waiting times in the GI/G/1 queue. The connection here
is not as direct as in the insurance risk model in that it is derived
from an analysis of ladder heights for transient random walks; see,
for example, Asmussen (1987, p.80), Feller (1971, p.396) and Pakes
(1975). Another common way in which this model occurs is through the
solution of a transient renewal equation. An example of this occurs in
branching processes, where we obtain a geometric-compound sum in the
analysis of the mean number of particles alive at a given time in an
age-dependent subcritical process; see Athreya and Ney (1972,
p.151). We refer to Feller (1971, chapter XI) for a discussion of
transient renewal theory. For further applications of subexponentiality
in transient renewal theory, we refer to Teugels (1975) and Embrechts
and Goldie (1982).

Throughout the paper, we assume that the $X_i$'s are nonnegative.

\bigskip


\section{Main results.}
If it exists, the hazard rate $h=F'/\oF$ yields
the representation of the distribution function $F$ as
$$
  \oF (t) = \oF(t_0)\exp\Bigl(-\int_{t_0}^t h(u)\d u\Bigr) \, .
$$
We write $\Id$ the identity function on $\RR$; for any positive real number
$r$, the function $\Id^r$ maps $t$ to $t^r$.
From the representation of $\oF$ with its hazard rate, we see that 
if $h\sim\alpha/\Id$ at infinity, then
$\oF$ is regularly varying with index $-\alpha$. If $\limt h(t)=\alpha$
then $\oF(t)=e^{-\alpha t(1+o(1))}$ has a tail behavior close to that of
an exponential distribution. Since we are interested in light subexponential
tails, it is natural to consider hazard rates such that
\setbox1=\vbox{\hsize=3.2in\par
  \noindent$h$ is regularly varying,\hfill\vskip3pt\break\noindent
  $\limt th(t)=+\infty \qquad \hbox{ and }\qquad \limt h(t)=0 \, .$
}
$$
  \lower8pt\box1\eqno{\equa{BasicHa}}
$$
In order to be not too close to the Pareto type distributions, we will 
strengthen this assumption by requiring that
$$
  \liminf_{t\to\infty} th(t)/\log t > 0 \, .
  \eqno{\equa{BasicHb}}
$$

This excludes distributions with tail $e^{-(\log t)^a}$ with $a< 2$,
but include those for which $a\geq 2$. It also includes the subexponential
Weibull distributions, or more generally, those with tail of the form $t^\beta
e^{-t^\alpha}$ with $\alpha$ positive and less than $1$.

As observed in Barbe and McCormick (2004, 2005), smoothness is a key 
requirement
to obtain asymptotic expansions. For our purposes, a good class of regularly
varying functions are the smoothly varying ones of given order, whose 
definition we now recall.

\bigskip

\Definition{ A function $h$ is smoothly varying of index $\alpha$ and order $m$
  if it is  ultimately $m$-times continuously differentiable and its $m$-th 
  derivative is regularly varying of index $\alpha-m$. 
}

\bigskip

Clearly, if the hazard rate is $m$ times differentiable, the tail function 
$\oF$ can be differentiated $m+1$ times.

\bigskip

The next notation we need to introduce pertains to the Laplace
characters. We write $\D$ the derivation operator; that is, if $g$ is 
differentiable, $\D g$ is its derivative. As is customary, we define $\D^0$ to
be the identity, and for any positive integer $i$ we define $\D^i$ by
induction as $\D \D^{i-1}$.

We write $\mu_{F,i}$ the $i$-th moment of $F$.

\bigskip

\Definition{ (Barbe and McCormick, 2004). Let $F$ be a distribution
  function having at least $m$ moments. Its Laplace character of order
  $m$ is the differential operator 
  $$
    \L_{F,m}=\sum_{0\leq i\leq m} {(-1)^i\over i!}\mu_{F,i}\D^i \, .
  $$
}

Laplace characters have useful algebraic properties which are
described in Barbe and McCormick (2004). In particular, consider the
ring $\RR_m[\, D\,]$ defined as the quotient ring of polynomials in
$\D$ modulo the ideal generated by $\D^{m+1}$. Laplace characters are
elements of this ring, and can be multiplied. It may be helpful to
think of a Laplace character as a formal Laplace transform $Ee^{-X\D}$
where $X$ has distribution $F$, expressed as a formal Taylor series in
$\D$, dropping all terms in $\D^{m+1},\D^{m+2}, \ldots$. Then, the
multiplication in the ring $\RR_m[\,\D\,]$ amounts to the usual
multiplication of Taylor series, dropping any term in
$\D^{m+1},\D^{m+2}, \ldots$ In particular, in $\RR_m[\,\D\,]$, we have
$\L_{H*K,m}=\L_{H,m}\L_{K,m}$.  In what follows, we always consider
Laplace characters of order $m$ as members of $\RR_m[\, D\,]$, and all
the operations on Laplace characters are in that quotient ring.

\medskip

The following theorem provides an asymptotic expansion for the tail of $G$.

\Theorem{\label{MainTheorem}
  Let $F$ be a distribution function whose hazard rate is smoothly varying 
  with negative index at least $-1$ and  positive order $m$. Assume 
  further that {\BasicHb} holds and that the moment generating function of
  $N$ is finite in a neighborhood of the origin. Then for any nonnegative
  integer $k$ at most $m$
  $$
    \oG=EN\L_{F^{\star (N-1)},k} \oF + o(h^k\oF) \, .
  $$%
  \vskip -8pt
}
\Remark It is shown in Barbe and McCormick (2005, Lemma 4.1.1) that
under the assumptions of Theorem {\MainTheorem} the asymptotic
equivalence $\oF^{(k)}\sim (-1)^kh^k\oF$ holds. Therefore, the
remainder term in the above formula could be written as $o(\oF^{(k)})$.

\bigskip


\section{Examples.} We illustrate the use of Theorem \MainTheorem, considering
the cases where $N$ has a Poisson and a geometric distribution.

\medskip

\noindent{\it Example 1.} Assume that $N$ has a Poisson distribution with
parameter $a$. Sums with a Poisson number of summands are commonly used
in insurance mathematics, modelling total claim size (see  Beirlant et al., 
1996, Embrechts et al., 1997, Willmot and Lin, 2000). The following expansion
is easily derived.

\Proposition{\label{PropPoisson}
  Let $F$ be a distribution function satisfying the assumptions of Theorem
  {\MainTheorem}. If $N$ has a Poisson distribution with parameter $a$, then
  $\oG=a\L_{G,m}\oF+o(h^m\oF)$. Moreover,
  $\L_{G,m}=e^{a(\L_{F,m}-\Id)}$.
}

\bigskip

\Proof Combine Theorem {\MainTheorem} and the proof of Corollary 4.4.2 in 
Barbe and McCormick (2004) to obtain the expansion $a\L_{G,m}\oF$. To obtain
the expression for $\L_{G,m}$, write, in the quotient ring,
$$
  EN\L_{F,m}^{N-1}
  = e^{-a}\sum_{n\geq 1} n {a^n\over n!} \L_{F,m}^{n-1}
  = ae^{a(\L_{F,m}-\Id)}\, .\eqno{\qed}
$$

\bigskip

\noindent The above formula is easily implemented with a computer algebra 
system. For example, the following {\tt Maple} code calculates 
$ae^{a(\L_{F,m}-\Id)}$.

\verbatim@
mu[0]:=1: 
LF:=sum('(-1)^j*mu[j]*D^j/j!','j'=0..m+1):
taylor(a*exp(a*(LF-1)),D=0,m+1);

@

Setting $m=3$ in the previous code yields the first four terms,
$$\eqalign{\qquad
  E(N\L_{F,4}^{N-1})=a\Id 
  &{}-a^2\mu_{F,1}\D + {a^2\over 2} (a\mu_{F,1}^2+\mu_{F,2})\D^2\cr
  &{}-{a^2\over 6}(a^2\mu_{F,1}^3+3a\mu_{F,1}\mu_{F,2}+\mu_{F,3})\D^3 
   \, .\cr
  }
$$

To give a very concrete example, assume that $\oF$ is the Weibull distribution
with parameter $1/3$, so that $\oF(t)=e^{-t^{1/3}}$. 
Define $e_r(t)=t^re^{-t^{1/3}}$. We obtain, after 
evaluation of $E(N\L_{F,4}^{N-1})$, and using a computer algebra 
package,
$$\displaylines{
  \oG= ae_0 
+ 2a^2e_{-2/3} 
+2a^2(20+a) e_{-4/3}
{}+ 4a^2(20+a)e_{-5/3}
\hfill\cr\noalign{\vskip 3pt}\hskip .7cm
{}+{4a^2(1680+60a+a^2)\over 3} e_{-2}
{}+8a^2\bigl(1680+60a+a^2)    e_{-7/3}
\hfill\cr\noalign{\vskip 2pt}\hskip .7cm
{}+{2 a^2(403200+9120a+140a^2+a^3)\over 3} e_{-8/3} +o(e_{-8/3})\, .
\hfill\cr}
$$
Perhaps the only remarkable feature of 
such 7 terms expansion is that it can be done.

\bigskip

\noindent{\it Example 2.} Motivated by applications to queueing theory (see
e.g., Cohen, 1972, or Bingham, Goldie and Teugels, 1987, p.387), consider the
case where $N$ has a geometric distribution with parameter $a$, that is
$N$ is a nonnegative integer $n$ with probability $(1-a)a^n$. Again, Theorem
{\MainTheorem} provides a compact expression of the asymptotic expansion
of $\oG$, and the issue is how to actually compute it. 

Any polynomial in $\D$ with nonvanishing constant term is invertible
in the quotient ring $\RR_m[\,\D\,]$. Therefore,
$$
  EN\L_{F^{\star (N-1)},m}=(1-a)\sum_{n\geq 1} a^n n \L_{F,m}^{n-1}
  = a(1-a)(\Id -a\L_{F,m})^{-2} \, .
$$
Consequently, the following result holds.

\Proposition{\label{PropGeometric}
  Let $F$ be a distribution function satistfying the assumptions of Theorem
  {\MainTheorem}. If $N$ has a geometric distribution with parameter $a$, then
  $\oG=a(1-a)(\Id-a\L_{F,m})^{-2}\oF+o(h^m\oF)$.
}

\bigskip

Setting $m=3$, we obtain, as in the previous example, with the help of a
computer algebra package, with $b=a/(1-a)$,
$$\eqalign{
  EN\L_{F,3}^{N-1}
  = b\,\Id 
  &{}- 2b^2\mu_{F,1}\D + b^2(\mu_{F,2}+3b\mu_{F,1}^2)\D^2 \cr
  \noalign{\vskip 2pt}
  &{}- {b^2\over 3} (12b^2\mu_{F,1}^3+9b\mu_{F,1}\mu_{F,2}+\mu_{F,3})\D^3 
  \, .\cr}
$$
 
For instance, when $F$ is the 
Weibull distribution with parameter $1/3$, the calculation of 
$EN\L_{F,4}^{N-1}$ yields the following 7 terms expansion
--- expressed solely with $a$, the formula contains alternating signs; it is 
numerically slightly more stable when expressed with $b=a/(1-a)$.

$$\displaylines{
  \oG= 
  b\, e_0 
  +4b^2\, e_{-2/3}
  +4b^2(20+3b)\, e_{-4/3}
  +8b^2(20+3b)\, e_{-5/3}
  \hfill\cr\noalign{\vskip 5pt}\hskip 20pt
  {}+32b^2(140+15b+b^2)\, e_{-2}
  {}+192b^2(140+15b+b^2)\, e_{-7/3}
  \hfill\cr\noalign{\vskip 5pt}\hskip 20pt
  {}+ 80b^2(6720+456b+28b^2+b^3)\, e_{-8/3} +o(e_{-8/3})\, .
  \hfill\cr}
$$

\bigskip


\section{Proof.}
When $m$ vanishes, Theorem {\MainTheorem} is due to Embrechts, Goldie and
Veraverbeke (1979, p.342). Therefore, we will prove it when $m$ is at 
least $1$.

It is convenient to introduce a pseudo-semi-norm on tails.
If $K$ is a distribution function, we write
$$
  {|\oK|}_F=\sup_{t\geq 0} (\oK/\oF)(t) \, ,
$$
with the convention $0/0=0$.  This generates balls $B(F,r)$ containing
all tails $\oK$ which are less than $r\oF$. We write $B(F)$ the union
of all these balls for all positive $r$.

We write $G_n$ for the $n$-fold convolution $F^{\star n}$. 

We start by recalling Kesten's global bound on tail function of 
self-convolutions of subexponential distributions; see Athreya and Ney (1972,
\S IV.4, Lemma 7). It asserts that
for any positive $\epsilon$ there exists a positive $A$ such that for all
positive integers $n$,
$$
  {|\oG_n|}_F\leq A(1+\epsilon)^n \, .
  \eqno{\equa{Kesten}}
$$

We also need a precise estimate of the order of magnitude of derivatives of
$\oF$. As noted in the Remark following Theorem \MainTheorem, Lemma 4.1.1 in 
Barbe and McCormick (2005) shows that for any nonnegative $k$ at most $m$,
$$
  \oF^{(k)}\sim (-1)^k h^k\oF \, .
  \eqno{\equa{DerivativeOrder}}
$$

Finally, we also need a basic representation of convolution in terms of
operators. For any distribution function $K$ with support in the nonnegative
half-line and any $\eta$ positive and less
than $1$, define the operator
$$
  \T_{K,\eta} f(t)=\int_0^{\eta t} f(t-x)\d K(x) \, .
$$
For any positive $c$ we also define the multiplication operator $\M_c$ acting 
on functions by
$$
  \M_c f(t)=f(t/c) \, .
$$
These two operators allow us to write a convolution in a way suitable for 
our analysis. Define the powers $\T_{K,\eta}^{\,n}$ by $\T_{K,\eta}^{\,0}=
\Id$ and $T_{K,\eta}^{\,n+1}=T_{K,\eta}\T_{K,\eta}^{\,n}$.
Using Proposition 5.1.1 in Barbe and McCormick (2004) 
inductively we obtain a representation for the distribution function,
valid on the nonnegative half line,
$$
\oG_n=\sum_{1\leq i\leq n} \T_{F,\eta}^{\,i-1}\T_{G_{n-i},1-\eta}\oF
  {}+\sum_{1\leq i\leq n}\T_{F,\eta}^{\,i-1}(\M_{1/\eta}\oF \M_{1/(1-\eta)}
  \oG_{n-i}) \, .
  \eqno{\equa{RepG}}
$$

Our first lemma is a simple moment bound.

\bigskip

\Lemma{\label{MomentBound}
  Let $i$ be a nonnegative integer, and let $\epsilon$ be a positive
  real number. There exists $t_1$ such that for
  any $t$ at least $t_1$ and any distribution function $K$ in $B(F)$,
  $$
    \int_t^\infty x^i \d K(x) \leq (1+\epsilon) {|\oK|}_F t^i\oF(t) \, .
  $$
}

\Proof For any nonnegative integer $i$, an integration by parts yields
$$
  \int_t^\infty x^i \d K(x)
  = t^i\oK(t)+i\int_t^\infty x^{i-1}\oK (x) \d x \, .
  \eqno{\equa{MomentBounda}}
$$
The right hand side of this equality is less than ${|\oK|}_F$ times the 
same expression with
$K$ replaced by $F$. Consequently, it suffices to prove the result
when $K$ is $F$. In that case, let $M$ be a positive real number so
that $\epsilon (M-i)\geq i$.  Since $\Id\, h$ tends to infinity at
infinity, $h$ is more than $M/\Id$ ultimately. For any $t$
large enough and any $x$ at least $t$,
$$
  {\oF(x)\over \oF(t)} = \exp\Bigl( -\int_t^xh(u)\d u\Bigr) 
  \leq \Bigl({t\over x}\Bigr)^M
  \, .
$$
This implies that the integral in the right hand side of
{\MomentBounda}, when $F$ is substituted for $K$, is at 
most $\epsilon t^i\oF(t)$.\hfill$\qed$

\bigskip

Our next lemma contains the main argument of the proof, namely that a
$\T_{K,\eta}$ operator is in some sense very close to a Laplace
character as far as tail behavior is concerned when applied to $\oF$
and its derivatives.

\bigskip

\Lemma{\label{TApprox}
  For any fixed integer $p$ at most $m$,
  $$
    \limt \sup_{K\in B(F)} 
    {|(\T_{K,\eta}-\L_{K,m-p})\oF^{(p)}|\over {|\oK|}_F h^m\oF } (t) = 0 \, .
  $$
}

\Proof The proof of Lemma 4.2.3 in Barbe and McCormick (2005) shows that
for any $\delta$ positive,
$$
  \Bigl|\int_0^{\delta/h(t)} \oF^{(p)}(t-x)\d K(x) - \L_{K,m-p}\oF^{(p)}\Bigr|
  \eqno{\equa{TApproxa}}
$$
is at most
$$\displaylines{
  \sum_{0\leq j\leq m-p} |\oF^{(p+j)}(t)|\int_{\delta/h(t)}^\infty x^j\d K(x)
  \hfill\equa{TApproxb}\cr\hfill
  {}+\int_0^{\delta/h(t)}\int_0^x {y^{m-p-1}\over (m-p-1)!} |\oF^{(m)}(t-x+y)
  - \oF^{(m)}(t)|\d y \d K(x)\, .
  \ \equa{TApproxc}\cr}
$$
Let $\epsilon$ be a positive number. Using Lemma \MomentBound\ and 
\DerivativeOrder, we see that for large $t$, the term \TApproxb\ is less than 
$$
  \oF(t) 2{|\oK|}_F \Bigl({\delta\over h(t)}\Bigr)^j \oF
  \Bigl( {\delta\over h(t)}\Bigr)
  \, .
$$
Since $\oF$ is rapidly varying, this is ultimately less than 
$\epsilon{|\oK|}_F h^m\oF$.

The proof of Lemma 4.2.3 in Barbe and McCormick (2005) shows that for
$\delta$ small enough, for any $t$ large enough and for any $K$ in
$B(F)$, the double integral \TApproxc\ is at most $\epsilon {|\oK|}_F
\mu_{F,m-p}h^m\oF$. Hence, we have shown that \TApproxa\ is at most
$\epsilon {|\oK|}_F(\mu_{F,m-p}+1)h^m\oF$ ultimately uniformly over
$B(F)$.

The proof of Lemma 4.2.4 in Barbe and McCormick (2005) shows that for any
positive $\delta$ and $\eta$, ultimately uniformly over $B(F)$,
$$
  \int_{\delta/h(t)}^{\eta t}|\oF^{(p)}(t-x)|\d K(x) 
  \leq \epsilon {|\oK|}_F h^m\oF(t) \, .
$$
This proves Lemma \TApprox.\hfill$\qed$

\bigskip

Lemma {\TApprox} yields the following estimate on an operator $\T$ composed
with a Laplace character applied to a derivative of $\oF$.

\Lemma{\label{TLApprox}
  The following uniform limit holds:
  $$
    \limt \sup_{\matrix{\ss K\in B(F)\cr\noalign{\vskip -3pt}
                        \ss H : \mu_{H,m-p}<\infty\cr}}
    { |\T_{K,\eta}\L_{H,m-p}\oF^{(p)}-\L_{K\star H,m-p}\oF^{(p)}|
      \over
      {|\oK|}_F\, h^m\oF\, \sum_{0\leq j\leq m-p}{\mu_{H,j}\over j!} 
    }
    = 0 \, .
  $$
}

\Proof Since $\T_{K,\eta}$ is linear and
$$
  \L_{H,m-p}\oF^{(p)}=\sum_{0\leq j\leq m-p} {(-1)^j\over j!} \mu_{H,j}
  \oF^{(p+j)}\, ,
$$
the result follows from Lemma {\TApprox} and Lemma 2.1.4 in Barbe and McCormick
(2004).\hfill$\qed$

\bigskip

The next two lemmas will take care of some remainder terms. The first one
asserts that terms of order $o(h^m\oF)$ remain so through the action of 
some $\T$ operators.

\bigskip

\Lemma{\label{LemmaRemainder}
  Let $q$ be a nonnegative integer and $\epsilon$ be a positive real number.
  There exist $t_2$, some positive $A$ and $\eta$, such that for
  any positive integer $i$,
  $$
    \T_{F,\eta}^{\,i}(h^q\oF) \leq A(1+\epsilon)^i h^q\oF
  $$
  on $[\,t_2,\infty)$.
}

\bigskip

\Proof Let $\epsilon$ be a positive real number.
Since $h$ is regularly varying with negative index, provided $\eta$ is
small enough, $h(t-x)\leq (1+\epsilon)h(t)$
for any $t$ large enough and any $x$ nonnegative and at most $\eta t$.

Therefore, for $t$ at least $t_2'$,
$$\eqalign{
  \T_{F,\eta}(h^q\oF)(t)
  &{}=\int_0^{\eta t} h^q\oF(t-x)\d F(x) \cr
  &{}\leq (1+\epsilon) h^q(t) \int_0^{\eta t} \oF(t-x)\d F(x) \cr
  &{}\leq (1+\epsilon) h^q(t)\oF^{\star 2}(t) \, . \cr
}
$$
By induction, it follows that
$$
  \T_{F,\eta}^{\,i} (h^q\oF)(t) 
  \leq (1+\epsilon)^i \bigl(h^q\oF^{\star (i+1)}\bigr)(t) \, .
$$
Using Kesten's bound, {\Kesten} above, this yields that
$\T_{F,\eta}^{\,i}(h^q\oF)$ is ultimately at most $A(1+\epsilon)^{2i} h^q\oF$,
finishing the proof since $\epsilon$ is arbitrary.\hfill$\qed$

\bigskip

Our penultimate lemma will be used to handle the terms involving the 
multiplication operators in {\RepG}.

\bigskip

\Lemma{\label{NeglectM}
  Let $\epsilon$ be a positive real number.
  There exists $t_3$ such that for any positive integers $i$ and $m$,
  $$
    |\M_{1/\eta}\oF\M_{1/(1-\eta)} \oG_i|
    \leq (1+\epsilon)^i h^{m+1}\oF
  $$
  on $[\, t_3,\infty)$.
}

\bigskip

\Proof Kesten's bound shows that
$$
  \bigl|\oF(t\eta)\oG_i\bigl( t(1-\eta)\bigr)\bigr|
  \leq \oF(t\eta) A(1+\epsilon)^i \oF\bigl( t(1-\eta)\bigr) \, .
$$
Arguing as in Lemma 4.2.1 in Barbe and McCormick (2005), 
$\oF(t\eta)\oF\bigl( t(1-\eta)\bigr)$ is $o\bigl( h^q\oF(t)\bigr)$ for any
positive $q$. This implies the result.\hfill$\qed$

\bigskip

Our last lemma is stated merely to avoid digression in the argument later on. 

\Lemma{\label{Boundmgf}
  Let $\epsilon$ be a positive number. There exists $A$ such that 
  for any positive integer $n$
  $$
    \sum_{0\leq j\leq m} {\mu_{G_n,j}\over j!} \leq A(1+\epsilon)^n \, .
  $$
}

\Proof The lemma folllows from Marcinkiewicz-Zygmund's inequality (see
Chow and Teicher, 1988, \S 10.3, Theorem 3), which implies that
$\mu_{G_n,j}\leq An^j$ for some constant $A$.\hfill$\qed$

\bigskip

We can now conclude the proof of Theorem \MainTheorem. Combining Lemmas
\LemmaRemainder\ and \NeglectM, there exists an interval $[\, t_3,\infty)$
on which for any $j$ and $k$ with $0\leq j\leq k\leq n$, any positive $i$
and $n$ with $i\leq n$,
$$
  |\T_{F,\eta}^{\,i-1}(\M_{1/\eta}\oF\M_{1/(1-\eta)} \oG_{n-i})|
  \leq A(1+\epsilon)^n h^{m+1} \oF \, .
$$
Representations {\RepG} yield, on $[\,t_3,\infty)$,
$$
  |\oG_n -\sum_{1\leq i\leq n} \T_{F,\eta}^{\,i-1}\T_{G_{n-i,1-\eta}}\oF|
  \leq An(1+\epsilon)^n h^{m+1} \oF  \, .
  \eqno{\equa{Maina}}
$$

Let $\epsilon$ be a positive real number, small enough so 
that $E (1+2\epsilon)^N$ is finite.
Let $\delta$ be a positive real number. Combining Lemmas  
{\TApprox}, {\TLApprox}, {\LemmaRemainder} and \Boundmgf, using also 
Kesten's bound, ultimately,
uniformly in $n$ and $i$ at most $n$,
\hfuzz=1pt
$$\displaylines{%
  |\T_{F,\eta}^{\,i-1}\T_{G_{n-i},1-\eta} \oF -\T_{F,\eta}^{\,i-2}
   \L_{G_{n-i+1},m} \oF|
  \hfill\cr\noalign{\vskip 3pt}
  {}\leq\T_{F,\eta}^{\,i-1}|(\T_{G_{n-i},1-\eta}-\L_{G_{n-i},m})\oF|
           + \T_{F,\eta}^{\,i-2} 
           |(\T_{F,\eta}\L_{G_{n-i},m}-\L_{G_{n-i+1,m}})\oF| 
  \hfill\cr\noalign{\vskip 3pt}
  {}\leq 2\delta (A^2+A)(1+\epsilon)^n h^m\oF \, .
  \hfill\equa{Mainb}\cr}
$$
\hfuzz=0pt
Using the same combination of lemmas, we also have, ultimately, uniformly
in $n$ and $j$ at most $n-1$,
$$
  |\T_{F,\eta}^{\,j}\L_{G_{n-j-1},m}\oF -\T_{F,\eta}^{\,j-1} \L_{G_{n-j},m}
   \oF|
  {}\leq A^2(1+\epsilon)^n\delta h^m\oF \, .
  \eqno{\equa{Mainc}}
$$
We take $A$ to be at least $1$, simply to ensure that $A^2$ is more than
$A$. Summing {\Mainc} for $j$ positive and less than $i$ and adding 
\Mainb, we obtain
$$
  |\T_{F,\eta}^{\,i-1}\T_{G_{n-i},1-\eta}\oF -\L_{G_{n-1},m}\oF|
  \leq 4A^2\delta h^m\oF i(1+\epsilon)^n 
$$
on some interval $[\,t_4,\infty)$. Summing these inequalities for $i$ positive
and at most $n$ and combining with {\Maina} yield
$$
  |\oG_n-n\L_{G_{n-1},m}\oF|
  \leq 10 A^2\delta n(n+1)  (1+\epsilon)^nh^m\oF \, .
$$
Since the moment generating function of $N$ 
is finite at $\log (1+2\epsilon)$ and $\delta$ is arbitrary, Theorem
{\MainTheorem} follows.\hfill$\qed$

\bigskip


\noindent{\bf References}
\medskip

{\leftskip=\parindent
 \parindent=-\parindent
 \par

S.\ Asmussen (1987). {\sl Applied Probability}, Wiley.

S.\ Asmussen (1997). {\sl Ruin Probabilities}, World Scientific.

K.B.\ Athreya, P.E.\ Ney (1972). {\sl Branching Processes.}
Springer.

Ph.\ Barbe, W.P.\ McCormick (2004). Asymptotic expansions for infinite
weighted convolutions of heavy tail distributions and applications,
{\tt http://www.arxiv.org/abs/math.PR/0412537}, submitted.

\hfuzz=3pt
Ph.\ Barbe, W.P.\ McCormick (2005). Asymptotic expansions for convolutions
of light tailed subexponential distributions, 
{\tt http://www. arxiv.org/abs/math.PR/0512141}, submitted

\hfuzz=0pt
J.\ Beirlant, J.L.\ Teugels, P.\ Vynckier (1996). {\sl Practical Analysis of 
Extreme Values}, Leuven University Press, Leuven, Belgium.

N.H.\ Bingham, C.M.\ Goldie, J.L.\ Teugels (1989) {\sl Regular Variation},
2nd ed., Cambridge.

D.B.H.\ Cline (1987). Convolutions of distributions with exponential and
subexponential tails, {\sl J.\ Austr.\ Math.\ Soc.\ (A)}, 43, 347--365.

Y.S.\ Chow, H.~Teicher (1988). {\sl Probability Theory, Independence, 
Interchangeability, Martingales}, 2nd ed., Springer.

J.W.\ Cohen (1972). On the tail of the stationary waiting-time distribution
and limit theorem for M/G/1 queue, {\sl Ann.\ Inst.\ H.\ Poincar\'e},{\sl B},
8, 255--263. 

P.\ Embrechts (1985). Subexponential distribution functions and their 
applications: a review, in {\sl Proceedings of the Seventh Conference on
Probability Theory, (Bra\c sov, 1982)}, 125--136, VNU Sci.\ Press, Utrecht, 
1985.

P.\ Embrechts, C.M.\ Goldie (1982). On convolution tails, {\sl Stoch.\ Proc.\
Appl.}, 13, 263--278.

P.\ Embrechts, C.M.\ Goldie, N.\ Veraverbeke (1979). Subexponentiality and
infinite divisibility, {\sl Z.\ Wahrsch.\ verw.\ Geb.}, 49, 335--347.

P.\ Embrechts, C.~Kl\"uppelberg, T.~Mikosch (1997). {\sl Modelling 
Extremal Events}, Springer.

W.\ Feller (1971). {\sl An Introduction to Probability and its Applications},
2nd ed., Wiley.

C.\ Goldie, C.\ Kl\"uppelberg (1998). Subexponential distributions, in 
{\sl A Practical Guide to Heavy Tails, Statistical Techniques and 
Applications}, R.\ Adler, R.\ Feldman and M.\ Taqqu eds., Birkh\"auser, 
435--460.

R.\ Gr\"ubel (1987). On subordinated distributions and generalized renewal
measures, {\sl Ann.\ Probab.}, 15, 394--415.

E.\ Omey, E.\ Willekens (1987). Second-order behaviour of distributions 
subordinate to a distribution with finite mean, {\sl Comm.\ Statist.\ Stoch.\
Models}, 3, 311-342.

A.G.\ Pakes (1975). On the tails of waiting-time distributions, {\sl J.\ Appl.\
Probab.}, 12, 555--564.

J.L.\ Teugels (1975). The class of subexponential distributions, {\sl Ann.\
Probab.}, 3, 1000--1011.

G.E.~Willmot, X.S.~Lin (2000). {\sl Lundberg Approximations for Compound
Distributions with Insurance Applications}, {\sl Lecture Notes in Statistics},
156, Springer.

}

\bigskip
\setbox1=\vbox{\halign{#\hfil&\hskip 40pt #\hfill\cr
  Ph.\ Barbe            & W.P.\ McCormick and C.\ Zhang\cr
  90 rue de Vaugirard   & Dept.\ of Statistics \cr
  75006 PARIS           & University of Georgia \cr
  FRANCE                & Athens, GA 30602 \cr
                        & USA \cr
                        & $\{$bill, czhang$\}$@stat.uga.edu \cr}}
\box1

\bye

\vfill\eject


{\bf This is not meant to stay in the paper. I wrote
it because I originally thought I would use it for the proof. I think it is
interesting though, because that gives the results of the other paper for
densities in a rather automatic paper.}

\bigskip

\section{Extra.} The class of distribution functions considered
is defined through representation of their tail in term of the hazard
function and smooth variation of the hazard rate. The purpose of this
subsection is to prove the following proposition, showing that this 
representation is stable by differentiation.

\Proposition{\label{Stability} 
  Let $F$ be a distribution function whose hazard rate $h$
  is smoothly varying of index $\alpha-1$ and order $m$, with $\alpha$ 
  positive and at most $1$. If $\alpha$ is $1$, assume further that
  $\limt th(t)=\infty$. For any nonnegative integer $p$ at most $m$
  there exists a function $h_p$ smoothly varyig of index $\alpha-1$ and 
  order $m-p$ such that
  $$
    \oF^{(p)}(t)=(-1)^p\oF(t_0)\exp\Bigl( -\int_{t_0}^t h_p(u)\d u\Bigr) \, .
  $$
  Moreover, $h_p\sim h$ at infinity.
}

\bigskip

The proof is based on the following lemma.

\Lemma{\label{StabilityLemma}
  Let $f$ be a function smoothly varying of nonpositive index $-\beta$ 
  and order $p$. If $\beta$ is less than $1$ or $\beta$ is $1$ and 
  $\limt tf(t)=\infty$, then $f+(\log f)'$ is smoothly varying of index
  $-\beta$ and order $p-1$. More precisely, for any $k$ at most $p$,
  the $k$-th derivative of $f+(\log f)'$ is asymptotically equivalent
  to that of $f$. 
}

\bigskip

\Proof Let $g$ be $f+(\log f)'$. If suffices to prove that $g^{(p-1)}$ is
regularly varying of index $-\beta-p+1$. Recall that for $k$ at least $1$,
the $k$-th derivative of the logarithm is $(-1)^{k-1}k!/\Id^k$. Therefore,
Fa\`a di Bruno's formula (see e.g.\ Roman, 1980) yields
$$
  (\log f)^{(p)}=\sum_{1\leq k\leq p} {(-1)^{k-1}\over f^k}
  \sum_{\matrix{\ss n_1+\cdots+n_k=p\cr\noalign{\vskip -3pt} 
                \ss n_1,\cdots,n_k\geq 1}}
  {p!\over n_1!\ldots n_k!}\, f^{(n_1)}\ldots f^{(n_k)} \, .
$$
Since $f$ is smoothly varying, 
$f^{(n_i)}\sim (-1)^{n_i}(\beta)_{n_i}f/\Id^{n_i}$. Therefore, if $n_1+\cdots
+n_k=p$ and the $n_i$'s are positive,
$$
  f^{(n_1)}\ldots f^{(n_k)} \sim f^k {(-1)^p\over \Id^p}(\beta)_{n_1}\ldots
  (\beta)_{n_k} \, .
$$
It follows that
$$\displaylines{\quad
  (\log f)^{(p)}
  \sim \Id^{-p} \sum_{1\leq k\leq p} {(-1)^{k-1}\over f^k}
  \hfill\cr\hfill
  \sum_{\matrix{\ss n_1+\cdots+n_k=p\cr\noalign{\vskip -3pt} 
                \ss n_1,\cdots,n_k\geq 1}}
  {p!\over n_1!\ldots n_k!}(\beta)_{n_1}\ldots (\beta)_{n_k} \, ,
  \quad\equa{EqStabilitya}\cr
  }
$$
with the usual convention that if the summation in the right hand side of
{\EqStabilitya} vanishes, {\EqStabilitya} means $(\log f)^{(p)}=o(\Id^{-p})$.
However, the right hand side in {\EqStabilitya} does not depend on $f$. In
the special case $f(x)=x^\alpha$, the left hand side of {\EqStabilitya} can
be explicitely calculated and is equal to $(-1)^p\beta p!/\Id^p$; hence the
right hand side of {\EqStabilitya} is $(-1)^p\beta p!/\Id^p$. This implies
$$
  g^{(p-1)}=f^{(p-1)}+(\log f)^{(p)} \sim (-1)^{p-1}(\beta)_{p-1}
  {f\over/\Id^{p-1}} \, .
$$
and proves that $g^{(p-1)}$ is regularly varying with 
index $\beta-p+1$.\hfill$\qed$

\bigskip

\noindent{\bf Proof of Proposition {\Stability}.} The proof is by induction,
the result being tautological when $k$ vanishes. If the representation holds
for a certain $p$, then
$$
  \oF^{(p+1)}=(-1)^{p+1} \exp\Bigl( -\int_0^1 h_p(u)+(\log h_p)'(u) \d u\Bigr)
  \, .
$$
We set $h_{p+1}=h_p+(\log h_p)'$ and apply Lemma {\StabilityLemma} to show
that $h_{p+1}$ is smoothly varying of same index as $h_p$ and of order 
$m-p-1$.$\qed$

\bigskip

\noindent{\bf References}
\medskip

{\leftskip=\parindent
 \parindent=-\parindent
 \par

S.\ Roman (1980). The formula of Fa\`a di Bruno, {\sl Amer.\ Math.\ Monthly},
87, 805--809.

}

\bye